\documentclass{article}

\usepackage{amsmath,amssymb,amsthm}
\usepackage{graphicx,psfrag}

\usepackage{pstricks}

\newtheorem{theorem}{Theorem}[section]

\newtheorem{lemma}[theorem]{Lemma}

\theoremstyle{remark}

\theoremstyle{definition}

\begin{document}

\title{Two-Dimensional Problems of Minimal Resistance
in a Medium of Positive Temperature\footnote{Research report
CM04/I-08, Dep. Mathematics, Univ. Aveiro, 2004.
Accepted to the Proceedings of the
Sixth Portuguese Conference on Automatic Control --
Controlo 2004, Faro, Portugal, June 7-9, 2004.}}

\author{Alexander Yu. Plakhov and Delfim F. M. Torres\\
        \texttt{\{plakhov,delfim\}@mat.ua.pt}}

\date{Department of Mathematics\\
      University of Aveiro\\
      3810-193 Aveiro, Portugal}

\maketitle


\begin{abstract}
We study the Newton-like problem of minimal resistance for a
two-dimensional body moving with constant velocity in a
homogeneous rarefied medium of moving particles. The distribution
of the particles over velocities is centrally symmetric. The
problem is solved analytically; the minimizers are shown to be of
four different types. Numerical results are obtained for the
physically significant case of gaussian circular distribution of
velocities, which corresponds to a homogeneous ideal gas of
positive temperature.
\end{abstract}


\smallskip



\section{Introduction}

In 1686, in his \emph{Principia Mathematica}, Newton considered
one of the oldest problems of optimal control. The problem
consists of finding the shape of a body, moving with constant
velocity in a homogeneous medium consisting of infinitesimal
particles, such that the total resistance of the medium to the
body would be minimal. Newton assumed that the collisions of the
particles with the body are absolutely elastic, the medium is very
rare, so that the particles do not mutually interact, and that the
particles are immovable, \textrm{i.e.}, there is no temperature
motion of particles. Problems of this kind may appear in
construction of high-speed and high-altitude flying vehicles, such
as missiles and artificial satellites.

Newton solved this minimization problem in the class of convex
axially symmetric bodies with the axis parallel to the velocity of
the body, of fixed length along the axis and fixed maximal cross
section orthogonal to the axis. Convexity implies that each
particle hits the body at most once, and this fact allows one to
write down the explicit formula for resistance.
An account of the problem originally considered by Newton can be found,
\textrm{e.g.}, in the book \cite{MR56:4953}.

Since the early days of the calculus of variations, many
modifications of Newton's problem of minimal resistance have been
studied in the literature \cite{JFM44.0455.02}. In 1990th
the interest to Newton's problem revived. Interesting results were
obtained when dealing with the minimization problem in various
classes of bodies obtained by withdrawing or relaxing the
conditions initially imposed by Newton:~ axial symmetry
\cite{BK,BFK,LP};~ convexity (still maintaining the
single shock assumption) \cite{CL,MR2003i:49030};~ single
shock condition \cite{P1,P3}.

On the other hand, the problem was studied under the more
realistic assumptions of presence of friction at the moment of
collision (so that the collisions are not absolutely elastic)
\cite{friction},\, and of mutual interaction of particles
\cite{W}. In the present paper, the case of temperature noise of
particles is considered. First, we obtain general formulas in the
$d$-dimensional case, $d \ge 2$, and then study the case $d = 2$
in more detail.

We suppose that a convex and axisymmetric body is moving in
${\mathbb R}^d$ with a constant velocity $V$, in a homogeneous
medium of moving particles; the distribution of particles over
velocities is the same at every point. In fact, it is more
convenient to assume that the body is immovable, and there is a
flux of particles falling upon it. This picture will be taken in
the sequel. The length $h$ of the body along the axis is fixed,
and the maximal cross section by a hyperplane orthogonal to the
axis is supposed to be a $(d-1)$-dimensional ball of radius 1.

In section~\ref{sec:2} the preliminary analysis of the problem in
the $d$-dimensional case is made. In section~\ref{sec:3}, the
minimization problem for $d = 2$ is solved. In
section~\ref{sec:4}, we consider the physically relevant special
case of gaussian circular distribution of velocities. Analytical
formulas for resistance are given, and results of numerical
simulations are presented.


\section{Calculation of Pressure and Resistance in the General Case}
\label{sec:2}

Consider a flux of infinitesimal particles in ${\mathbb R}^d$,\,
$d \ge 2$. The density of flux and the distribution of particles
over velocities are the same at each point, the distribution being
given by a density function $\rho(v)$. The pressure of the flux at
a regular point $x \in
\partial \mathcal{B}$ of the boundary of a convex body
$\mathcal{B} \subset {\mathbb R}^d$ equals $\pi(n_x)$,
where $n_x$ is the outer normal to $\partial \mathcal{B}$ at $x$,
\begin{equation}\label{pi}
\pi(n) = - \int (v | n)_{\!-}^{\,\ 2} \rho(v)\, dv \cdot n;
\end{equation}
here $x_- := \min \{ x,\, 0 \}$ and $(\cdot | \cdot)$ means scalar
product; and resistance of the body to the flux equals
\begin{equation}\label{R}
R(\mathcal{B}) = \int_{\partial \mathcal{B}} \pi(n_x)\, d{\cal H}^{d-1}(x),
\end{equation}
where ${\cal H}^{d-1}$ means $(d-1)$-dimensional Hausdorff measure.
The formulas (\ref{pi}) and (\ref{R}) are obvious modifications of
the corresponding formulas from \cite{BG}.

In what follows, we shall assume that

(i) the function $\rho$ is spherically symmetric with the center
$-V e_d$, where $V > 0$ and $e_d$ is the $d$th coordinate vector,
i.e. $\rho(v) = \sigma(|v + V e_d|)$;\, besides

(ii) the function $\sigma$ is continuous and monotone decreasing, and
$$
\int_0^\infty r^2 \sigma(r)\, r^{d-1} dr < \infty.
$$

The body $\mathcal{B}$ is supposed to be convex, open, and symmetric with
respect to the $d$th coordinate axis. By translation along this
axis and by subsequent homothety, it can be reduced to the form
\begin{equation}
\label{B}
\mathcal{B} =
\{ (x\text{'}, x_d) : |x\text{'}| < 1,\ f_-(|x\text{'}|) < x_d
< -f_+(|x\text{'}|)\}\, ,
\end{equation}
$\mathcal{B} \subset {\mathbb R}^d$,
where $x\text{'} = (x_1, \ldots, x_{d-1})$,\, $f_+$ and $f_-$ are
convex negative non-decreasing functions defined on $[0,\, 1)$. We
shall imagine that the $d$th coordinate axis is directed
vertically upwards. The height of the body is
\begin{equation}\label{h=}
h = -f_+(0) - f_-(0).
\end{equation}
We shall suppose that $h$ is constant.

At a point $x_+ = (x\text{'}, -f_+(|x\text{'}|))$ of the upper
part of the boundary $\partial \mathcal{B}$, the outer normal equals
\begin{equation}\label{nx+}
n_{x_+} = \frac{1}{\sqrt{f_+'(|x\text{'}|)^{\,2} + 1}}\
\left(f_+'(|x\text{'}|) \frac{x\text{'}}{|x\text{'}|},\ 1 \right),
\end{equation}
and due to the property (ii) of $\rho$, from (\ref{pi}) one finds
that pressure of the flux at this point equals $\pi(n_{x_+}) =
-p_+ \left(f_+'(|x\text{'}|) \right) \cdot n_{x_+}$, where
\begin{equation}\label{p+}
p_+(u) := \left| \pi \left(\left( \frac{u}{\sqrt{u^2+1}},\, 0,
\ldots, 0,\, \frac{1}{\sqrt{u^2+1}} \right)\right)\right|.
\end{equation}

Similarly, the outer normal to $\partial \mathcal{B}$ at a point $x_- =
(x\text{'}, f_-(|x\text{'}|))$ of the lower part of $\partial \mathcal{B}$
equals $\pi(n_{x_-}) = p_- \left(f_-'(|x\text{'}|) \right) \cdot
n_{x_-}$, where
\begin{equation}
\label{p-}
p_-(u) :=
-\left| \pi \left(\left( \frac{u}{\sqrt{u^2+1}},\, 0,
\ldots, 0,\, -\frac{1}{\sqrt{u^2+1}} \right)\right)\right|.
\end{equation}

From (\ref{p+}), (\ref{p-}), and (\ref{pi}) one obtains
\begin{equation}\label{formula}
p_\pm(u) = \pm \int \frac{(v_1 u \pm v_d)\!_-^{\,\ 2}}{1 + u^2}\
\rho(v)\, dv.
\end{equation}

Using the assumptions (i) and (ii), one can show that the
functions $p_+$ and $p_-$ satisfy the following properties:

(a) $p_\pm \in C^1[0,\, +\infty)$;

(b) there exist $\lim_{u\to +\infty} p_\pm(u)$;

(c) $p_\pm'(0) = \lim_{u\to +\infty} p_\pm'(u) = 0$;

(d)  $p_+'(u) < p_-'(u) \le 0$.

(e) There exist $\bar{u}_\pm > 0$ such that\, $p_\pm'$ is strictly
monotone decreasing on $[0,\, \bar{u}_\pm]$, and strictly monotone
increasing on $[\bar{u}_\pm,\, +\infty)$.

The proof of these properties is not very difficult, but rather
technical and long, and will be presented elsewhere.

Let us calculate the resistance $R(\mathcal{B})$, using the formula
(\ref{R}). The integral in the right hand side of (\ref{R}) is the sum of
two integrals corresponding to the upper and the lower parts of
$\partial \mathcal{B}$. Change the variable in each of these integrals,
substituting $d{\cal H}^{d-1}(x)$ for $\sqrt{ f_+'(|x\text{'}|)^2 +
1}\, dx\text{'}$ and for $\sqrt{ f_-'(|x\text{'}|)^2 + 1}\,
dx\text{'}$, and substituting $\pi(n_{x})$ for
$\pi(n_{x_+})$ and for $\pi(n_{x_-})$, respectively. One obtains
\begin{multline*}
R(\mathcal{B})
= \int_{|x\text{'}|<1} p_+ (f_+'(|x\text{'}|)) \cdot
\left(f_+'(|x\text{'}|) \frac{x\text{'}}{|x\text{'}|},\ 1
\right)\, dx\text{'} \\
+ \int_{|x\text{'}|<1} p_-
(f_-'(|x\text{'}|)) \cdot \left( -f_-'(|x\text{'}|)
\frac{x\text{'}}{|x\text{'}|},\ 1 \right)\, dx\text{'},
\end{multline*}
and as a result of integrating, using the fact that the functions
$p_\pm(f_\pm'(|x\text{'}|))$ and $\pm f_\pm'(|x\text{'}|)
\frac{x\text{'}}{|x\text{'}|}$ are symmetric and antisymmetric,
respectively, with respect to $x\text{'}$, one gets
\begin{equation*}
R(\mathcal{B}) = \left( \int_{|x\text{'}|<1} p_+ (f_+'(|x\text{'}|))\,
dx\text{'}
+ \int_{|x\text{'}|<1} p_- (f_-'(|x\text{'}|))\,
dx\text{'} \right)\ e_d \, .
\end{equation*}
Therefore
\begin{equation}\label{RRR}
R(\mathcal{B}) = a_{d-1}\, \left( {\cal R}_+(f_+) + {\cal R}_-(f_-) \right) \cdot
e_d,
\end{equation}
where $a_{d-1}$ is the volume of a unit ball in ${\mathbb R}^{d-1}$, and
\begin{equation}\label{rpm}
{\cal R}_\pm (f) = \int_0^1 \, p_\pm (f'(t))\, dt^{d-1}.
\end{equation}

From the formulas (\ref{RRR}), (\ref{rpm}), and (\ref{h=}) it is
seen that the problem of finding the functions $f_+$,\, $f_-$,
minimizing $R(\mathcal{B})$, can be solved in two steps. First, find the
values
$$
\inf_{f\in{\cal M}(h_\pm)} {\cal R}_\pm(f)
$$
in the class ${\cal M}(h_\pm)$ of negative convex functions $f$ such
that $f(0) = -h_\pm$.\ Second, minimize the sum $\inf_{{\cal M}(h_+)}
{\cal R}_+(f) + \inf_{{\cal M}(h_-)} {\cal R}_-(f)$ provided that $h_+ + h_-
= h$.

Let us fix the sign "+" or "$-$", and introduce shorthand notation
$p_\pm = p, \ \ \ h_\pm = h, \ \ \ {\cal R}_\pm = {\cal R}$.

The following auxiliary lemma is a consequence
of the Pontryagin Maximum Principle \cite{MR29:3316b}.

\begin{lemma}\label{l3}
Let $\lambda > 0$,\, $f_h \in {\cal M}(h)$,\, $f_h(1) = 0$, and let for
any point $t$ of differentiability of $f_h$, $u = f_h'(t)$ be a
solution of the problem
\begin{equation}\label{min}
t^{d-2}\, p(u) + \lambda\, u \to \min.
\end{equation}
Then $f_h$ is a unique solution of the minimization problem
\begin{equation}\label{minim}
\inf_{f\in{\cal M}(h)} {\cal R}(f), \ \ \ \ \ {\cal R} (f) = \int_0^1 \, p
(f'(t))\, dt^{d-1}.
\end{equation}
\end{lemma}

\begin{proof} (\textrm{cf.} \cite{T}).
Indeed, for any $f \in {\cal M}(h)$,\, $f \ne f_h$ one has
$$
t^{d-2}\, p(f'(t)) + \lambda\, f'(t) \ge t^{d-2}\, p(f_h'(t)) +
\lambda\, f_h'(t).
$$
Moreover, this inequality is strict on a set of $t$ of positive
measure. Integrating both sides of this inequality over $t \in
[0,\, 1]$, one gets
\begin{multline*}
\frac{1}{d-1}\, \int_0^1 \, p (f'(t))\, dt^{d-1} + \lambda\, (f(1) -
f(0)) \\
> \frac{1}{d-1}\, \int_0^1 \, p (f_h'(t))\, dt^{d-1} +
\lambda\, (f_h(1) - f_h(0)),
\end{multline*}
and using that $f(1) \le f_h(1) = 0$ and $f(0) = f_h(0) = -h$,
one obtains that ${\cal R} (f) > {\cal R} (f_h)$.
\end{proof}


\section{Solution of the Problem in the Two-Dimensional Case}
\label{sec:3}

From now on, we shall assume that $d = 2$.

Using the property (e) of the function $p$, it is easy to prove
that the problem
$$
\frac{p(0) - p(u)}{u} \to \max
$$
has a unique solution; denote it by $u^0$.

Also, denote
$$
B := \frac{p(0) - p(u^0)}{u^0} = - p'(u^0).
$$
There may appear three different cases:

(a) As $\lambda > B$, the unique solution of (\ref{min}) is $u = 0$.

(b) As $\lambda = B$, there are two solutions:\, $u = 0$ and $u =
u^0$.

(c) As $\lambda < B$, the solution $\tilde u$ is unique, besides
$\tilde u > u^0$, and $p'(\tilde u) = -\lambda$.

Consider these cases separately.

(a) $\lambda > B$. The unique solution of (\ref{min}) is $u = 0$,
hence, according to lemma \ref{l3}, $f \equiv 0$ is the solution
of (\ref{minim}) for $h = 0$, and ${\cal R}(f) = p(0)$.

(b) $\lambda = B$. There are two solutions:~ $u = 0$ and $u =
u^0$, hence a function $f_h \in {\cal M}(h)$, whose derivative
takes the values 0 and $u^0$, is the solution of (\ref{minim}). By
virtue of convexity of $f_h$, there exists $t_0 \in [0,\ 1]$ such
that $f_h'(t) = 0$ as $t \in [0,\ t_0]$,\, and $f_h'(t) = u^0$ as
$t \in [t_0,\, 1]$.\, Thus,
\begin{equation}\label{xt}
f_h(t) = \left\{ \begin{array}{lc} -h & \text{ as } t \in [0,\ t_0]\\
-h + u^0 \cdot (t - t_0)  & \text{ as } t \in [t_0,\, 1].
\end{array} \right.
\end{equation}
Taking into account that $f_h(1) = 0$, one concludes that $h \le
u^0$ and $t_0 = 1 - h/ u^0$. Further, one has
\begin{equation*}
{\cal R}(f_h) = \int_0^{t_0} p(0)\, dt + \int_{t_0}^1 p(u^0)\, dt
= p(0) + \frac{h}{u^0}\, (p(u^0) - p(0)).
\end{equation*}
Using the definition of $B$, one gets
$$
{\cal R}(f_h) = p(0) - B\, h \, .
$$

(c) $\lambda < B$. There is a unique solution $\tilde u$ of
(\ref{min}), hence the function
\begin{equation}\label{tx}
f_h(t) = -h + \tilde u t
\end{equation}
solves (\ref{minim}), where $h = \tilde u$. Here the resistance
equals
$$
{\cal R}(f_h) = p(h) \, .
$$

Define the function
$$
\bar{p}(u) = \left\{ \begin{array}{ll} p(0) - B\, u & \ \ \text{
if } \ 0 \le u \le u^0,\\
p(u) & \ \ \text{ if } \ u \ge u^0.
\end{array} \right.
$$
Thus, one comes to

\begin{lemma}\label{l4}
The solution $f_h$ of (\ref{minim}) is given by(\ref{xt}), if $h <
u^0$, and by (\ref{tx}), if $h \ge u^0$; moreover, ${\cal R}(f_h)
= \bar p(h)$.
\end{lemma}

Reverting to the subscripts "+" and "$-$" and using lemma
\ref{l4}, one concludes that the problem of finding
$$
\inf_{h_+ + h_- = h} \left( {\cal R}_+(f_{h_+}) + {\cal
R}_-(f_{h_-}) \right) =: \mathrm R(h)
$$
amounts to the problem
\begin{equation}\label{minp}
\min_{0 \le z \le h} p_h(z), \ \ \ \text{where} \ \ \ p_h(z) =
\bar p_+(z) + \bar p_-(h - z).
\end{equation}

The functions $\bar p_\pm(u)$ are continuously differentiable on
$[0,\, +\infty)$, and their derivatives are monotone increasing,
hence $p_h'(z)$,\, $0 \le z \le h$\, is also monotone increasing.

Using the property (d), one concludes that $B_+ > B_-$, hence
there exists a unique value $u_* > u_+^0$ such that $\bar
p_+'(u_*) = -B_-$. Consider four cases:

1) $h < u_+^0$;

2) $u_+^0 \le h \le u_*$;

3) $u_* < h < u_* + u_-^0$;

4) $h \ge u_* + u_-^0$.

In the cases 1) and 2) one has $p_h'(z) \le p_h'(h) = \bar p_+'(h)
+ B_- \le 0$ as $0 \le z \le h$, hence $z = h$ is the solution of
(\ref{minp}). Therefore, the optimal value of $h_-$ is zero, and
$f_{h_-}(t) \equiv 0$.

1) $h < u_+^0$. One has $h_+ = h < u_+^0$, hence $f_{h_+}$ is
given by (\ref{xt}), with $t_0 = 1 - h/u_+^0$.\, The optimal body
is a trapezium, with tangent of slope of its lateral sides equal
to $u_+^0$. The minimal resistance equals $\mathrm R(h) = p_+(0) -
B_+\, h + p_-(0)$.

2) $u_+^0 \le h \le u_*$. One has $f_{h_+}(t) = -h + h\, t$, hence
the optimal body is an isosceles triangle. Here $\mathrm R(h) =
p_+(h) + p_-(0)$.

In the cases 3) and 4)\, one has $\bar p_+'(h) > -B_- > -B_+ =
\bar p_+'(u_+^0)$, hence $h > u_+^0$.\, Further, one has $p_h'(h)
= \bar p_+'(h) - B_- > 0$; on the other hand, $p_h'(u_+^0) = \bar
p_+'(u_+^0) - \bar p_-'(h - u_+^0) \le B_- - B_+ < 0$. It follows
that the minimum of $p_h$ is achieved at an interior point of
$[u_+^0,\, h]$, so $u_+^0 < h_+ < h$,\, and $f_{h_+}(t) = (t -
1)\, h_+$.

3) $u_* < h < u_* + u_-^0$. Denoting $\tilde h = \max \{ 0,\, h -
u_-^0 \}$, one has $\tilde h < u_*$, hence
$$
p_h'(\tilde h) = \bar p_+'(\tilde h) - \bar p_-'(h - \tilde h)
\le \bar p_+'(\tilde h) + B_- < 0 \, ,
$$
therefore the minimum of $p_h$ is reached at an interior point of
$[\tilde h,\, h]$. Thus, $0 < h_- < h - \tilde h \le u_-$, and
$f_{h_-}(t) = -h_-$ if  $t \in [0,\ 1 - h_-/u_-^0]$; $f_{h_-}(t) =
-h_- + u_-^0\, (t - 1 + h_-/u_-^0)$ if $t \in [1 - h_-/u_-^0,\
1]$. The optimal body here is the union of a triangle and a
trapezium turned over. The tangent of slope of lateral sides of
the trapezium equals $-u_-^0$. The minimal resistance equals
$\mathrm R(h) = p_+(u_*) + p_-(0) - B_-\, (h - u_*)$.

4) $h \ge u_* + u_-^0$. One has $p_h'(h - u_-^0) = \bar p_+'(h -
u_-^0) + B_- \ge 0$,\, hence the minimum of $p_h$ is reached at a
point of $[u_+^0,\, h - u_-^0)$.\, Thus, $h_- > u_-^0$, and
$f_{h_-}(t) = t\, h_-$.\, The optimal body is a union of two
isosceles triangles with common base and of heights $h_+$ and
$h_-$ defined from the relations $h_+ + h_- = h$,\, $p_+'(h_+) =
p_-'(h_-)$,\, $h_+ \ge u_+^0$,\, $h_- \ge u_-^0$. The minimal
resistance here equals $\mathrm R(h) = p_+(h_+) + p_-(h_-)$.


\section{Gaussian Distribution of Velocities -- Exact Solutions}
\label{sec:4}

Suppose that the density $\rho$ is gaussian circular, with the
mean $-Ve_d$ and variance 1, i.e.
\begin{equation}\label{gauss}
\rho_V(v) = \frac{1}{2\pi}\, e^{-\frac 12 |v + Ve_d|^2}.
\end{equation}
Here, the value $V$ is allowed to vary, so we shall denote the
pressure functions by $p_\pm(u, V)$ instead of $p_\pm(u)$. Fixing
the sign "+" and passing to polar coordinates $v = (-r
\sin\varphi, -r \cos\varphi)$ in the formula (\ref{formula}), one
obtains
\begin{equation}
\label{rhod=211}
p_+(u, V) = \int\!\!\int \frac{r^2 (\cos\varphi +
u\sin\varphi)_{\!+}^{\,\ 2}}{1 + u^2}\, \rho_+ (r, \varphi, V)\, r
dr d\varphi \, ,
\end{equation}
where $x_+ := \max \{ x,\, 0 \}$, and
$\rho_+(r, \varphi, V)$ is the gaussian density
(\ref{gauss}) written in the introduced polar coordinates,
\begin{equation}\label{rhod=222}
\rho_+(r, \varphi, V) = \frac{1}{2\pi}\, e^{-\frac 12 (r^2 - 2Vr
\cos\varphi + V^2)}.
\end{equation}
Next, fixing the sign "$-$" and introducing polar coordinates in a
slightly different way, $v = (-r \sin\varphi, r \cos\varphi)$, one
obtains
\begin{equation}
\label{rhod=233}
p_-(u, V) = -\int\!\!\int \frac{r^2 (\cos\varphi +
u\sin\varphi)_{\!+}^{\,\ 2}}{1 + u^2}\, \rho_- (r, \varphi, V)\, r
dr d\varphi \, .
\end{equation}
Here $\rho_-(r, \varphi, V)$ is the same density (\ref{gauss})
written in these coordinates,
\begin{equation}\label{rhod=244}
\rho_-(r, \varphi, V) = \frac{1}{2\pi}\, e^{-\frac 12 (r^2 + 2rV
\cos\varphi + V^2)}.
\end{equation}
Combining the formulas (\ref{rhod=211}), (\ref{rhod=222}),
(\ref{rhod=233}), and (\ref{rhod=244}), one comes to the more
general expression
\begin{equation}
\label{rhod=2} p_\pm(u, V) = \pm \frac{e^{-V^2/ 2}}{2\pi}
\int\!\!\hspace{-8mm}\int\limits_{\cos\varphi + u\sin\varphi > 0}
\hspace{-9mm}\left\{\frac{\left(\cos\varphi + u\sin\varphi\right)^2}{1 + u^2}
\, e^{-\frac 12 r^2 \pm 2rV \cos\varphi} r^3 \right\} dr
d\varphi \, .
\end{equation}

Passing to the iterated integral and integrating over $r$, one
obtains
\begin{equation*}
p_\pm(u, V) = \pm \frac{e^{-V^2/2}}{\pi}
\int\limits_{\cos\varphi + u\sin\varphi > 0}
\frac{(\cos\varphi + u\sin\varphi)^2}{1 + u^2}\ l(\pm V
\cos\varphi)\, d\varphi \, ,
\end{equation*}
where
\begin{equation}
\label{asympt1}
l(z) = 1 + \frac{z^2}{2}
+ \frac{\sqrt\pi}{2\sqrt 2} e^{z^2/ 2}
\left(3z + z^3\right) \left( 1 + \text{erf} \left( \frac{z}{\sqrt
2} \right) \right).
\end{equation}
Changing the variable $\tau = \varphi - \arcsin(u/ \sqrt{1+u^2})$,
one finally comes to
\begin{equation}\label{asympt2}
p_\pm(u, V) = \pm \frac{e^{-V^2/ 2}}{\pi} \int_{-\pi/2}^{\pi/2}
\cos^2 \tau\ l( \pm z(\tau,u,V) )\, d\tau,
\end{equation}
where $z(\tau,u,V) = V\, (\cos\tau - u\sin\tau)/ \sqrt{1 + u^2}$.

Using (\ref{asympt2}), one comes to the asymptotic formula
$$
p_\pm(u, V) = \pm \frac 12 + \sqrt{2/\pi}\, \frac{V}{\sqrt{1 +
u^2}} + O(V^2), \ \ V \to 0+.
$$
On the other hand, using (\ref{asympt1}), in the limit $V \to
\infty$ one obtains the asymptotic formula for pressure
$$
p_+(u, V) = \frac {V^2}{1 + u^2}\, (1 + o(1)), \ \ p_-(u, V) =
o(V^2);
$$
here the main term coincides, up to the factor $V^2$, with the
pressure in Newton (i.e. zero-temperature) case.

\begin{figure}
\begin{center}
\psline(3.25,2.2)(3,1.8)(3.5,1.8)(3.25,2.2) 
\psline(3.25,3.78)(3,3.28)(3.5,3.28)(3.25,3.78) 
\psline(3,3.28)(3.1,3.13)(3.4,3.13)(3.5,3.28) 
\psline(3.25,7.5)(3,6.5)(3.5,6.5)(3.25,7.5) 
\psline(3,6.5)(3.25,6)(3.5,6.5) 
\psline(3.1,1.15)(3,1.04)(3.5,1.04)(3.4,1.15)(3.1,1.15) 
\psfrag{v}{$V$}
\psfrag{h}{$h$}
\psfrag{0}{$0$}
\psfrag{1}{$1$}
\psfrag{2}{$2$}
\psfrag{3}{$3$}
\psfrag{4}{$4$}
\psfrag{5}{$5$}
\psfrag{10}{$10$}
\psfrag{15}{$15$}
\psfrag{20}{$20$}
\psfrag{25}{$25$}
\includegraphics[width=8cm,height=10cm]{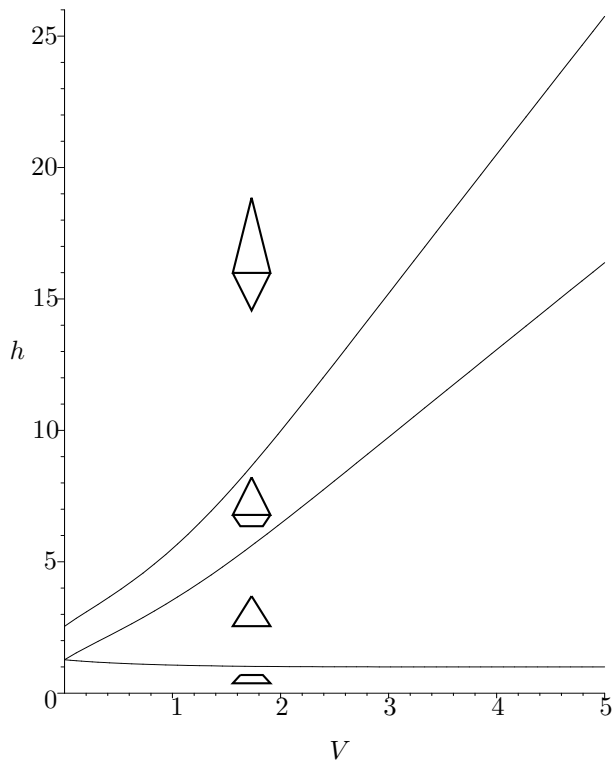}
\end{center}
\caption{Height $h$ versus velocity $V$ of the body}
\label{fig:vh}
\end{figure}

\begin{figure}
\begin{center}
\psfrag{h}{$h$}
\psfrag{RR}{$\hspace*{-7mm}\tilde{\mathrm R}(V,h)$}
\psfrag{0}{$0$}
\psfrag{0.5}{$0.5$}
\psfrag{1}{$1$}
\psfrag{1.5}{$1.5$}
\psfrag{2}{$2$}
\psfrag{2.5}{$2.5$}
\psfrag{3}{$3$}
\psfrag{3.5}{$3.5$}
\psfrag{2}{$\hspace*{-1mm}2$}
\psfrag{4}{$\hspace*{-1mm}4$}
\psfrag{6}{$\hspace*{-1mm}6$}
\psfrag{8}{$\hspace*{-1mm}8$}
\psfrag{10}{$\hspace*{-1mm}10$}
\psfrag{12}{$\hspace*{-1mm}12$}
\psfrag{14}{$\hspace*{-1mm}14$}
\psfrag{16}{$\hspace*{-1mm}16$}
\psfrag{v=0.1}{$V = 0.1$}
\psfrag{v=0.2}{$V = 0.2$}
\psfrag{v=0.5}{$V = 0.5$}
\psfrag{v=1}{$V = 1$}
\psfrag{v=infinity}{$V = \infty$}
\includegraphics[width=6.5cm,height=10cm]{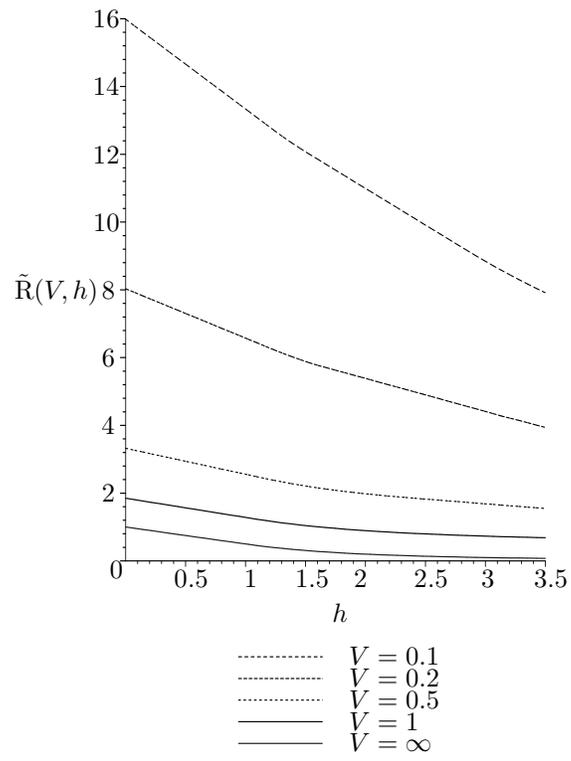}
\end{center}
\caption{Reduced resistance $\tilde{\mathrm R}(V,h)$
versus height $h$ of the body}
\label{fig:RRh1}
\end{figure}

\begin{figure}
\begin{center}
\psfrag{h}{$h$}
\psfrag{RR}{$\hspace*{-7mm}\tilde{\mathrm R}(V,h)$}
\psfrag{0}{$0$}
\psfrag{1}{$\hspace*{-1mm}1$}
\psfrag{2}{$2$}
\psfrag{3}{$3$}
\psfrag{4}{$4$}
\psfrag{5}{$5$}
\psfrag{6}{$6$}
\psfrag{0.2}{$\hspace*{-1mm}0.2$}
\psfrag{0.4}{$\hspace*{-1mm}0.4$}
\psfrag{0.6}{$\hspace*{-1mm}0.6$}
\psfrag{0.8}{$\hspace*{-1mm}0.8$}
\psfrag{1.2}{$\hspace*{-1mm}1.2$}
\psfrag{v=2}{$V = 2$}
\psfrag{v=3}{$V = 3$}
\psfrag{v=5}{$V = 5$}
\psfrag{v=infinity}{$V = \infty$}
\includegraphics[width=6.5cm,height=10cm]{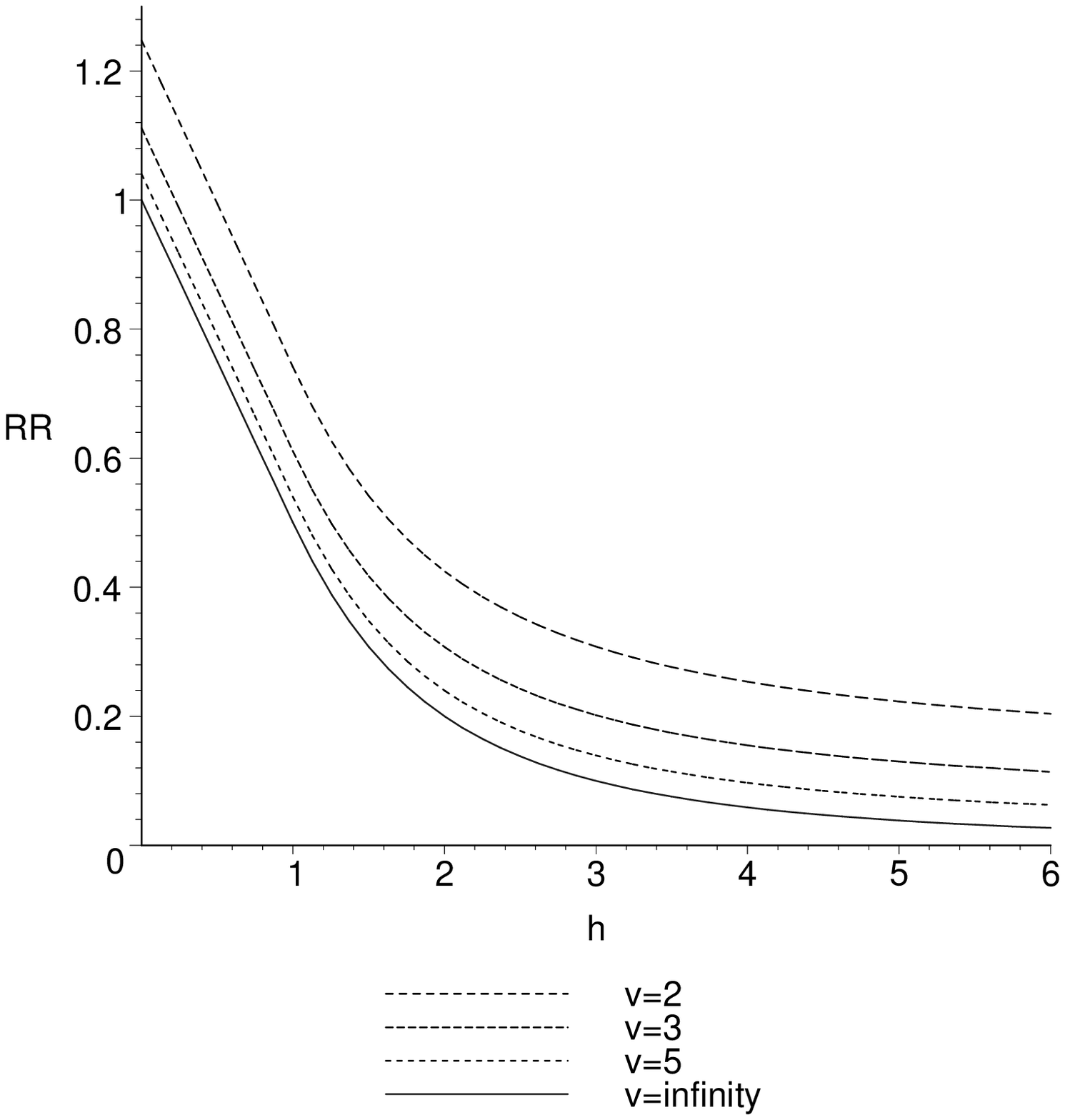}
\end{center}
\caption{Reduced resistance $\tilde{\mathrm R}(V,h)$
versus height $h$ of the body}
\label{fig:RRh2}
\end{figure}

Numerical simulations were done using Maple.
The results are presented on figures~\ref{fig:vh},
\ref{fig:RRh1} and \ref{fig:RRh2}.

Graphs of the three functions, shown on figure 1, separate the
$V$-$h$ plane into four regions, which correspond to the four
possible solutions. The lower function tends to 1 as $V \to
\infty$. As $V = 0$, the lower, the middle, and the upper
functions take the values $a$,\, $a$,\, and $2a$, respectively,
where $a = \sqrt{(1 + \sqrt 5)/ 2} \approx 1.27$.

Let $\mathrm R(V, h)$ be the minimal value of resistance at given
values $V$ and $h$; define reduced resistance by $\tilde{\mathrm
R}(V, h) = \mathrm R(V, h)/ V^2$. The graphs of $\tilde{\mathrm
R}(V, h)$ versus $h$ are shown on figures \ref{fig:RRh1} and
\ref{fig:RRh2} for different values of $V$.
As $V \to \infty$,\, $\tilde{\mathrm R}(V, h)$ tends to
$$
\tilde{\mathrm R}(\infty, h) = \left\{
\begin{array}{ll} 1 - h/2 & \text{ if } h \le 1\\
1/(1 + h^2) & \text{ if } h \ge 1.
\end{array}
\right.
$$
As $V \to 0$,\, $\tilde{\mathrm R} (V, h)$ goes to infinity,
besides
$$
\sqrt{2/\pi} \cdot \lim_{v\to\infty} (V\, \tilde{\mathrm R} (V,
h)) = \left\{
\begin{array}{ll} 2 - h/a^5 & \text{ if } h \le 2a\\
4/(4 + h^2) & \text{ if } h \ge 2a.
\end{array}
\right.
$$


\section{Conclusions and Final Comments}

In this paper, we treat one of the earliest problems in optimal
control: Newton's problem of minimal resistance. We have obtained
a full analytical description of the case when a two-dimensional
body moves through a rarefied medium of infinitesimal particles,
whose velocities are distributed according to the gaussian law.
From the physical viewpoint, such a medium is just a homogeneous
gas of positive temperature, while the case of immovable particles
is related to a zero-temperature gas.

The analytical formulas of the two-dimensional problem considered
in this work are quite complicated, even for Maple, but using the
current computational power, numerical simulations are made within
a few days of calculations. The study of the tree-dimensional case
would be a natural next step. However, calculation of the
resistance functions and the whole analysis are more involved.
This question is under development, and will be addressed
elsewhere.


\section*{Acknowledgements}

This work was partially supported by the R\&D
Unit CEOC of the University of Aveiro, and the
FCT/FEDER Project POCTI/MAT/41683/2001.


\end{document}